\documentclass[11pt,english]{article}
\usepackage[T1]{fontenc}
\usepackage[latin9]{inputenc}
\usepackage{textcomp}
\usepackage{amsthm}
\usepackage{amsmath}
\usepackage{amssymb}

\makeatletter
\theoremstyle{plain}
\newtheorem{thm}{Theorem}
  \theoremstyle{plain}
  \newtheorem{cor}[thm]{Corollary}

\usepackage{mathrsfs}\usepackage{amsthm}\usepackage{amscd}

\usepackage{hyperref}
\usepackage{graphics}
\usepackage{a4wide}
\@ifundefined{definecolor}
 {\usepackage{color}}{}

\newcounter{EQNR}
\renewcommand{\contentsname}{Contents\\
{\footnotesize\normalfont(A table of contents should normally not
be included)}}

\makeatother

\usepackage{babel}

\begin{document}

\title{Two extensions of Thurston's spectral theorem for surface diffeomorphisms}

\author{Anders Karlsson%
\thanks{Supported in part by Institut Mittag-Leffler (Djursholm, Sweden) and
the Swiss NSF grant 200021 132528/1.%
}}
\maketitle
\begin{abstract}
Thurston obtained a classification of individual surface homeomorphisms
via the dynamics of the corresponding mapping class elements on Teichmüller
space. In this paper we present certain extended versions of this,
first, to random products of homeomorphisms and second, to holomorphic
self-maps of Teichmüller spaces. 
\end{abstract}

\section{Introduction}

Let $M$ be an oriented closed surface of genus $g\geq2$. Let $\mathcal{S}$
denote the isotopy classes of simple closed curves on $M$ not isotopically
trivial. For a Riemannian metric $\rho$ on $M$ and aclosed curve
$\beta$, let $l_{\rho}(\beta)$ be the infimum of the length of curves
isotopic to $\beta.$

In a seminal preprint from 1976 \cite{T88}, Thurston classified surface
diffeomorphisms as being isotopic either to a periodic one, or else
reducible or pseudo-Anosov. A version of this was obtained earlier
in a series of papers by Nielsen \cite{N,N44}. Using the theory of
foliations of surfaces, Thurston showed the following consequence;
the proof is worked out in exposé 11 of \cite[Théorème Spectral]{FLP79}: 
\begin{thm}
(\cite[Theorem 5]{T88}) \label{thm:Thurston}For any diffeomorphism
$f$ of $M$, there is a finite set $1\leq\lambda_{1}<\lambda_{2}<...<\lambda_{K}$
of algebraic integers such that for any $\alpha\in\mathcal{S}$ there
is a $\lambda_{i}$ such that for any Riemannian metric $\rho$, \[
\lim_{n\rightarrow\infty}l_{\rho}(f^{n}\alpha)^{1/n}=\lambda_{i}.\]
 The map $f$ is isotopic to a pseudo-Anosov map iff $K=1$ and $\lambda_{1}>1.$
\end{thm}
This statement is analogous to the dynamical behaviour of linear maps
$A$ of finite dimensional vector spaces: the limits $\lim_{n\rightarrow\infty}\left\Vert A^{n}v\right\Vert ^{1/n}$
exist for every vector $v$, as is immediate from the Jordan normal
form. In this note we will obtain a few extensions of Theorem \ref{thm:Thurston}.
First:
\begin{thm}
\label{thm:randomNT}For any integrable ergodic cocycle $f_{n}=g_{n}g_{n-1}...g_{1}$
of orientation-preserving homeomorphisms of $M$, there are almost
surely a constant $\lambda\geq1$ and a (random) measured foliation
$\mu$ such that for any $\alpha\in\mathcal{S}$ with $i(\mu,\alpha)>0$
and Riemannian metric $\rho$,\[
\lim_{n\rightarrow\infty}l_{\rho}(f_{n}\alpha)^{1/n}=\lambda.\]

\end{thm}
Let $\mathcal{T}$ be the Teichmüller space of $M$ and for $x\in\mathcal{T}$,
the corresponding hyperbolic length of $\alpha\in\mathcal{S}$ is
denoted by $l_{x}(\alpha)$. The previous theorem is a direct consequence
of the more precise statement:
\begin{thm}
\label{thm:randomNTfine}In the setting of the previous theorem, we
have almost surely that for any $x\in\mathcal{T}$, there is an explicit
constant $C(\mu,x)>0$ such that for any $\epsilon>0$ there is a
number $N$ for which \[
C(\mu,x)i(\mu,\alpha)(\lambda-\epsilon)^{n}\leq l_{x}(f_{n}\alpha)\leq l_{x}(\alpha)(\lambda+\epsilon)^{n}\]
holds for any $\alpha\in\mathcal{S}$ and any $n>N.$ 
\end{thm}
Notice in particular the uniformity in $\alpha$ and the appearance
of $i(\mu,\alpha).$ 

The mapping class group $MCG(M)$ is the group of isotopy classes
of orientation-preserving homeomorphisms (or diffeomorphisms) of $M$:\[
MCG(M)=Homeo^{+}(M)/Homeo_{0}(M),\]
which acts by automorphisms of $\mathcal{T}(M)$. Thus every random
product of homeomorphism gives rise to a random product of mapping
classes and acts on $\mathcal{T}$. Kaimanovich-Masur \cite{KM96}
studied random walks on $MCG$. They proved that if the support of
the random walk measure generates a non-elementary subgroup, then
a.e. trajectory converges to points in the set of uniquely ergodic
foliations of $\mathcal{PMF}$. Taking this into account, we can deduce:
\begin{cor}
\label{cor:RW}Let $f_{n}=g_{n}g_{n-1}...g_{1}$ be a product of random
homeomorphisms where $g_{i}$ are chosen independently and distributed
with a probability measure of finite first moment and that generates
a subgroup containing two independent pseudo-Anosov maps. Then there
is a number $\lambda>1$ such that a.s. for any $\alpha\in\mathcal{S}$
and metric $\rho$ \[
\lim_{n\rightarrow\infty}l_{\rho}(f_{n}\alpha)^{1/n}=\lambda.\]

\end{cor}
This can be viewed as analogous to a well-known theorem of Furstenberg
and to Oseledets' multiplicative ergodic theorem for random products
of matrices. Information about these results and other references
to the vast literature on random dynamical systems can be found in
\cite{A98}.

Given a complex structure $x$ on $M$, the extremal length of an
isotopy class of curve $\alpha$ is\[
Ext_{x}(\alpha)=\sup_{\rho\in\left[x\right]}\frac{l_{\rho}(\alpha)^{2}}{Area(\rho)}\]
where the supremum is taken over all metrics in the conformal class
of $x$. Miyachi \cite{Mi08} noted that a normalized version of the
extremal length function, denoted $E_{P}$, extends continuously to
the whole Gardiner-Masur compactification, see section 3 below for
more details. A boundary point $P$ is called \emph{uniquely ergodic
}if $E_{P}(\beta)>0$ for all $\beta\in\mathcal{S}$.

By a theorem of Royden from 1971, later extended by Earle-Kra to surfaces
$M$ with punctures, the mapping class group, with some lower genus
exceptions, is isomorphic to the complex automorphism group of $\mathcal{T\mathrm{(M)}}$.
Here is a statement about more general holomorphic self-maps, thus
in a sense providing a certain extension of Theorem \ref{thm:Thurston}:
\begin{thm}
\label{thm:holomorphicNT}Let $f:\mathcal{T}\rightarrow\mathcal{T}$
be a holomorphic map and $x\in\mathcal{T}$. Then there is a number
$\lambda\geq1$ and a point $P$ in the Gardiner-Masur compactification
such that for all $n\geq1$ and any curve $\beta\in\mathcal{S}$,
\[
Ext_{f^{n}x}(\beta)\geq\left(\inf_{\alpha}\frac{Ext_{x}^{1/2}(\alpha)}{E_{P}(\alpha)}\right)^{2}E_{P}(\beta)^{2}\lambda^{n}\]
and, provided that $E_{P}(\beta)>0$,\[
Ext_{f^{n}x}(\beta)^{1/n}\rightarrow\lambda.\]

\end{thm}
The following can be seen as a weak generalization of the Nielsen-Thurston
classification of mapping classes to general holomorphic self-maps
of Teichmüller spaces:
\begin{thm}
\label{thm:DW}Let $f:\mathcal{T}\rightarrow\mathcal{T}$ be a holomorphic
map. Then either every orbit in $\mathcal{T}$ is bounded, or every
orbit leaves every compact set and there are associated points $P$
in the Gardiner-Masur boundary . If P is uniquely ergodic, then it
is unique and every orbit converges to this point in either compactification
and for some $\lambda\geq1$, any $\alpha\in\mathcal{S}$ and any
$x\in\mathcal{T}$$ $\[
Ext_{f^{n}x}(\alpha)^{1/n}\rightarrow\lambda\:\:\:\textrm{and }\:\:\inf_{\alpha}\frac{Ext_{f(x)}^{1/2}(\alpha)}{E_{P}(\alpha)}\geq\lambda\inf_{\alpha}\frac{Ext_{x}^{1/2}(\alpha)}{E_{P}(\alpha)}.\]

\end{thm}
This classifies $f$ as either having bounded orbits or else having
certain associated boundary points $P$ ({}``reducible'' vs pseudo-Anosov
depending on the nature of $P$). This is reminiscent of the Wolff-Denjoy
theorem in complex dynamics that, together with a theorem of Fatou,
classifies holomorphic self-maps of the unit disk. The $P$s can informally
be regarded as {}``virtual fixed points at infinity''. Examples
of important holomorphic self-maps beyond the automorphisms are the
Thurston skinning map in three-dimensional topology and the Thurston
pull-back maps in complex dynamics, see \cite{Mc90,Pi01,S11} for
more details. 

\textbf{Method of proof:}

For the ergodic statements, the techniques of Ledrappier and myself
\cite{KaL11} are being employed, although extended to asymmetric
metrics, such as Thurston's Lipschitz metric. Together with a comparison
with Teichmüller's metric, and, crucially, the new insights about
horofunctions by Cormac Walsh \cite{W11}, this leads to the main
result Theorem \ref{thm:randomNTfine}. The proof of the corollary
relies in addition on \cite{KM96}.

For Theorems \ref{thm:holomorphicNT} and \ref{thm:DW}, the starting
points is the well-known fact that the Techmüller metric coincides
with the Kobayashi metric, which implies that holomorphic maps are
$1$-Lipschitz in this metric (this is the only way in which the holomorphy
is being used), after that the proof uses recent results on horofunctions
\cite{Mi08,LS10,Mi11} and \cite{Ka01}.

\textbf{Further comments:}

- While Nielsen studied lifts of the diffeomorphisms to the hyperbolic
disk and boundary circle, Thurston instead compactified $\mathcal{T}$
in a natural way, by adding $\mathcal{PMF}$, and applied Brouwer's
fixed point theorem. Bers gave an alternative approach using more
classical theory, in particular the Teichmüller metric \cite{B78}.
For the spectral theorem (Theorem \ref{thm:Thurston}) Thurston used
foliation theory, see \cite{FLP79} and also \cite[§14.5]{FMa12}.
Our approach here instead uses two metrics on $\mathcal{T}$: Thurston's
and Teichmüller's. In this context it might be useful to point out
that the recent paper \cite{LPST11} contians a study parallel to
Bers' paper but for the Thurston metric.

- Following the work of Kaimanovich-Masur \cite{M95,KM96} it is natural
to wonder about so-called ray approximation, see \cite{K00} where
this problem is mentioned. The first result in this direction is a
paper of Duchin \cite{D05} which established ray approximation in
the Teichmüller metric for times when the ray visits the thick part.
Recently there is an elegant argument by Tiozzo \cite{Ti12} removing
the thick part requirement, this was first obtained for finitely supported
measures and then later (after the present paper was completed) under
finite first moment. This provides another approach to Corollary \ref{cor:RW}.
Ray approximation for Weil-Petersson geodesics in the more general
ergodic setting follows from the arguments in a paper by Margulis
and myself, see \cite{KaL06}. A recent analysis of the harmonic measure
can be found in \cite{Ga09}.

- During the last few years there has been great progress in proving
that almost every element in $MCG$ is pseudo-Anosov thanks to the
works of Maher, Rivin, Kowalski, and Lubotzky-Meiri, see e.g. \cite{Ma11}.
It is perhaps interesting to compare this with Corollary \ref{cor:RW}
above which states that random walk trajectories eventually look pseudo-Anosov
from the perspective of Theorem \ref{thm:Thurston}.

- Looking at Theorem \ref{thm:Thurston} or Oseledets' theorem in
the linear case, one might hope for a more precise statement than
Theorem \ref{thm:randomNT}. There are however some obvious obstructions
to keep in mind. The exponent $\lambda$ can in general take any real
value $\geq1$ since one can take an appropriately asymmetric measure
on a single pseudo-Anosov and its inverse. Moreover, for predicting
which $\alpha$ grows exponentially in length, one should bear in
mind the case of a coboundary $g(x)=f(\omega)f(T\omega)^{-1}$ and
that the cocycle might consist of reducible maps with respect to one
and the same curve system.

- Thurston used iteration on $\mathcal{T}$ in order to prove invariant
structures (\cite{Mc90}), for example the existence of a hyperbolic
three-manifolds fibering over the circle. In their paper \cite{KM96},
Kaimanovich-Masur deduced that subgroups of $MCG$ cannot be higher
rank lattices (see also \cite{FM98}). One could hope that the present
results also will find application.

- Our approach might be relevant for the study of random walks on
$Out(F_{n})$; the problem of extending the Kaimanovich-Masur theory
is a well-known research problem.

- As already indicated, one could moreover hope that our weak Wolff-Denjoy
analog could shed some light on the dynamics of rational maps in complex
dynamics, the so-called Thurston obstruction theorem also leads to
considering the iteration of a holomorphic self-map such as Thurston's
pull-back maps. In this latter context the fixed point corresponds
to the conformal map being combinatorially equivalent to a rational
map. This is not always the case, and the conditions for checking
this, i.e. verifying that the holomorphic self-map has bounded orbit
is difficult in practice. The curves $\alpha$ for which $E_{P}(\alpha)=0$
might be tightly related to the Thurston obstruction.

- Does every holomorphic self-map of Teichmüller space with bounded
orbit have a fixed point? A positive answer to this question would
be an extension of Nielsen realization and also strengthen the analog
of Theorem \ref{thm:DW} with the Wolff-Denjoy theorem.

\textbf{Acknowledgements:}

This work was done during a wonderful term spent at the Institute
Mittag-Leffler. Special thanks are due to the staff for their very
efficient and kind assistance. For this paper I was inspired by the
recent progress in the understanding of horofunctions by Walsh, Miyachi,
and Liu-Su. I got the idea of generalizing Theorem \ref{thm:Thurston}
to random products while preparing a presentation of \cite{T88} at
a reading seminar in Geneva and while keeping in mind similiar questions
that Margulis suggested to me years ago. In addition, I am grateful
for various discussions with Moon Duchin, Sebastian Hensel, Vadim
Kaimanovich, Athanase Papadopoulos, John Smillie, and Giulio Tiozzo.

\section{Proof of Theorem \ref{thm:randomNTfine}}

Let $(\Omega,p)$ be a standard Borel space with $p(\Omega)=1$ and
$T:\Omega\rightarrow\Omega$ an ergodic measure preserving transformation.
Assume that $g:\Omega\rightarrow Homeo^{+}(M)\twoheadrightarrow MCG(M)$
is a measurable map and let \[
Z_{n}(\omega)=g(\omega)g(T\omega)...g(T^{n-1}\omega)\]
which is called an \emph{ergodic cocycle}. Notice here that we have
shifted the order, so that in the terminology of the theorem, $f_{n}^{.}=Z_{n}^{-1}$
and the $g_{i}=\left(g(T^{i-1}\omega)\right)^{-1}$. A \emph{random
walk} on $MCG$ is the special case when the increments $g(T^{i}\omega)$
are assumed to be independent and identically distributed (this is
precisely the case when $(\Omega,p)$ is a product space and $T$
the shift). 

For $x\in\mathcal{T}$ denote by $l_{x}(\alpha)$ the minimal length
in its isotopy class in the hyperbolic metric $x$ (or more precisely
the class of isometric metrics all mutually isotopic). Let us recall
Thurston's asymmetric metric (\cite{T86}),\[
L(x,y)=\log\sup_{\alpha\in\mathcal{S}}\frac{l_{y}(\alpha)}{l_{x}(\alpha)}.\]

It is easy to see that $L$ verifies the triangle inequality, and
it also true that it seperates points although this is non-trivial.
Therefore $L$ satisfies all the axioms for a metric except the symmetry,
which indeed fails except in very special cases of surfaces with symmetries.
The triangle inequality reads\[
L(x,z)\leq L(x,y)+L(y,z).\]
The metric is clearly invariant: $L(gx,gy)=L(x,y)$ for $g\in MCG.$
The topology induced by $L$ coincides with the usual one, see \cite{PT07}.

Fix a base point $x_{0}\in\mathcal{T}.$ We will assume that \[
\int_{\Omega}L(g(\omega)x_{0},x_{0})+L(x_{0},g(\omega)x_{0})dp(\omega)<\infty,\]
in which case we refer to $f_{n}$ as an \emph{integrable ergodic
cocycle}, or informally as a \emph{random product} of mapping class
elements.

One has the following subadditivity property:\[
L(f_{n+m}(\omega)x_{0},x_{0})\leq L(f_{n}(\omega)f_{m}(T^{n}\omega)x_{0},f_{n}(\omega)x_{0}))+L(f_{n}(\omega)x_{0},x_{0})\]

\[
=L(f_{m}(T^{n}\omega)x_{0},x_{0})+L(f_{n}(\omega)x_{0},x_{0}).\]
From the subadditive ergodic theorem of Kingman one then knows that
for a.e. $\omega$ the following limit exists (the finite non-negative
value being independent of $\omega$ by the ergodicity assumption):\[
l:=\lim_{n\rightarrow\infty}\frac{1}{n}L(f_{n}(\omega)x_{0},x_{0}).\]
We now introduce concepts from the work of Walsh \cite{W11} which
gave a crucial inspiration for the present paper. It will be important
to consider functions $h$ in the so-called horofunction compactification
of $\mathcal{T}$, that is, for $\mu\in PMF$ \[
h_{\mu}(x)=\log\sup_{\alpha}\frac{i(\mu,\alpha)}{l_{x}(\alpha)}-\log\sup_{\beta}\frac{i(\mu,\beta)}{l_{x_{0}}(\beta)},\]
(notice that it is well-defined for projective equivalence classes
of measured foliations $\mu$) and for $x_{n}\rightarrow\mu$ in the
Thurston compactification one has\[
h_{\mu}(x)=\lim_{n\rightarrow\infty}L(x,x_{n})-L(x_{0},x_{n}).\]
These functions together with $h_{z}(x)=L(x,z)-L(x_{0},z)$ constitute
a compact space $H$ homeomorphic to Thurston\textasciiacute{}s compactification
as proved by Walsh. Notice that $ $\[
h_{z}(x)=L(x,z)-L(x_{0},z)\leq L(x,x_{0})\]
by the triangle inequality.

For $f\in MCG$ and $h$ as above let $F(g,h)=-h(g^{-1}x_{0}).$ We
note the following cocycle property:\[
F(g_{1},g_{2}h)+F(g_{2},h)=-(g_{2}\cdot h)(g_{1}^{-1}x_{0})-h(g_{2}^{-1}x_{0})\]
\[
=-h(g_{2}^{-1}g_{1}^{-1}x_{0})+h(g_{2}^{-1}x_{0})-h(g_{2}^{-1}x_{0})=F(g_{1}g_{2},h).\]

Note that moreover\[
L(gx_{0},x_{0})=-L(g^{-1}x_{0},g^{-1}x_{0})+L(x_{0}g^{-1}x_{0})=\max_{h\in H}F(g,h),\]
as follows from the triangle inequality. Following \cite{KaL06},
see also \cite{KaL11}, one introduces the skew product system $\bar{T}:\Omega\times H\rightarrow\Omega\times H$
by $\bar{T}(\omega,h)=(T\omega,g^{-\text{1}}(\omega)h)$ and checks
that with $\bar{F}(\omega,h):=F(g(\omega)^{-1},h)$ one has that\[
F(f_{n}(\omega)^{-1},h)=\sum_{i=0}^{n-1}\bar{F}(\bar{T}^{i}(\omega,h)).\]
Moreover, we have $\left|F(g^{.1}(\omega),h)\right|\leq\max\left\{ L(x_{0},g(\omega)x_{0}),L(g(\omega)x_{0},x_{0})\right\} $
so $F$ is integrable (notice here again a small difference due to
the asymmetric nature of $L$). From this point on, the proof runs
as in the references mentioned, that is, a special measure is constructed
which accounts for the drift as well as projecting to the original
one. Birkhoff's ergodic theorem is applied and a measurable section
is taken. We get that for a.e. $\omega$ there is an $h=h^{\omega}$
such that \begin{equation}
\lim_{n\rightarrow\infty}-\frac{1}{n}h(Z_{n}x_{0})=l.\label{eq:horo}\end{equation}
If $l>0$, which is the nontrivial case, it is immediate that $h$
is a boundary point, and hence by Walsh's theorem cited above, it
is of the form $h_{\mu}$ for some $\mu\in PMF.$ This means that
for every $\epsilon>0$ there is an $N$ such that for all $n>N$
one has \[
\log\sup_{\alpha}\frac{i(\mu,\alpha)}{l_{Z_{n}x_{0}}(\alpha)}-\log\sup_{\beta}\frac{i(\mu,\beta)}{l_{x_{0}}(\beta)}\leq-(l-\epsilon)n.\]
Letting $C_{\mu}^{-1}=\sup\frac{i(\mu,\beta)}{l_{x_{0}}(\beta)}$
we then obtain\[
\sup_{\alpha}\frac{i(\mu,\alpha)}{l_{Z_{n}x_{0}}(\alpha)}\leq C_{\mu}^{-1}e^{-(l-\epsilon)n},\]
which leads to that for every $\alpha$ we have\[
l_{Z_{n}x_{0}}(\alpha)\geq C_{\mu}i(\mu,\alpha)e^{(l-\epsilon)n}.\]

On the other hand, one knows (see for example \cite{CR07,LPST10})
that for any two points $x$ and $y$ in the thick part of Teichmüller
space up to an additive constant $L(x,y)\asymp L(y,x)$, thanks to
the symmetry of the Teichmüller metric. Therefore, since the orbit
$Z_{n}x_{0}$ stays in the thick part (as it does not move in the
moduli space), we have from the subadditive ergodic theorem mentioned
above that for all sufficiently large $n$,\[
L(x_{0},Z_{n}x_{0})=\log\sup_{\alpha}\frac{l_{Z_{n}x_{0}}(\alpha)}{l_{x_{0}}(\alpha)}\leq(l+\epsilon)n.\]
In view of these two inequalities we have for every $\alpha\in\mathcal{S}$
with $i(\mu,\alpha)>0$ that\[
l_{Z_{n}x_{0}}(\alpha)^{1/n}\rightarrow\lambda:=e^{l}.\]
Note that $l_{x_{0}}(f_{n}\alpha)=l_{Z_{n}x_{0}}(\alpha)$ and that
since $M$ is compact from the point of view of exponential growth
every Riemannian metric is equivalent, so we can replace $x_{0}$
with $\rho$. This finishes the proof of Theorems \ref{thm:randomNT}
and \ref{thm:randomNTfine}. We also remark is that if we do not assume
ergodicity, by ergodic decomposition, there is no change other than
that the {}``Lyapunov exponent'' $\lambda$ is now random and not
necessarily constant. A final comment is that $M$ could be allowed
to have punctures, except possibly for restricting to Riemannian metrics
equivalent to hyperbolic metrics.

\subsection*{Additional arguments in the random walk case}

Assume that $\nu$ is a probability measure on $MCG$ of finite first
moment (here meaning that the cocycle is integrable) and whose support
generates a non-elementary subgroup. A non-elementary subgroup is
by definition a subgroup which do not fix any finite set in $\mathcal{PMF}$.
It is shown in \cite{KM96} that this is equivalent to that there
are at least two pseudo-Anosovs with disjoint fixed point sets. In
particular it implies that the subgroup is non-amenable. Kaimanovich
and Masur \cite{KM96} established that random walks a.s. converge
to points in the set of uniquely ergodic foliations of $\mathcal{PMF}$.
Moreover, they show that the measure that this convergence gives rise
to, the harmonic or hitting measure, is the unique $\nu-$stationary
measure. In my work with Ledrappier we construct special $\nu-$stationary
measures in a general setting, see \cite[Theorem 18]{KaL11}. By uniqueness,
this measure is hence the same as the hitting one and by the remark
towards the end of the proof of \cite[Theorem 18]{KaL11} we get the
same conclusion (\ref{eq:horo}) above. Now the random foliations
$\mu$ are uniquely ergodic, so $i(\mu,\alpha)>0$ for any $\alpha\in\mathcal{S}$.
Finally, $\lambda>1$ because $l>0$ which in turn is a consequence
of entropy theory, including a well-known inequality of Guivarc'h,
and that a non-elementary subgroup is non-amenable as well as the
fact that $\mathcal{T}$ has at most exponential growth as established
in \cite{KM96}, see for example \cite{K00} for this type of arguments.
This proves Corollary \ref{cor:RW}.

\section{Proof of Theorem \ref{thm:holomorphicNT}}

Recall Kerckhoff's formula for the Teichmüller distance:\[
d(x,y)=\frac{1}{2}\log\sup_{\alpha\in\mathcal{S}}\frac{Ext_{x}(\alpha)}{Ext_{y}(\alpha)},\]
where $Ext$ denotes extremal length and defined in the introduction.
Royden showed that this metric coincides with the Kobayashi metric
associated to the complex structure of $\mathcal{T}.$ This implies
that holomorophic maps are 1-Lipschitz in this metric, see \cite{E01}
for an update on finer contraction properties. Gardiner and Masur
introduced in \cite{GM91} a compactification analogous to the Thurston
compactification but using extremal lengths instead of hyperbolic
lengths.

Fix a point $x_{0}\in\mathcal{T}.$ Let\[
E_{x}(\alpha)=\frac{Ext_{x}(\alpha)^{1/2}}{K_{x}^{1/2}},\]
where $K_{x}$ is the quasi-conformal dilation of the Teichmüller
map from $x_{0}$ to $x$. Miyachi \cite{Mi08} noted that $E$ extends
continuously to a function defined on the Gardiner-Masur compactification
$\overline{\mathcal{T}}^{GM}$of $\mathcal{T}.$

Let $f:\mathcal{T}\rightarrow\mathcal{T}$ be a holomorphic self-map
of Teichmüller space. Define \[
l=\lim_{n\rightarrow\infty}\frac{1}{n}d(f^{n}x_{0},x_{0})\]
which exists by the 1-Lipschitz property and subadditivity. Clearly,
$0\leq l<\infty$. Moreover, for any point $P\in\overline{\mathcal{T}}^{GM}$
define following Liu and Su\[
h_{P}(x)=\log\sup_{\beta}\frac{E_{P}(\beta)}{Ext_{x}(\beta)^{1/2}}-\log\sup_{\alpha}\frac{E_{P}(\alpha)}{Ext_{x_{0}}(\alpha)^{1/2}}.\]
We set for short \[
Q(P)=\sup_{\alpha}\frac{E_{P}(\alpha)}{Ext_{x_{0}}(\alpha)^{1/2}}.\]
Given a sequence $\epsilon_{i}\searrow0$ we set $b_{i}(n)=d(f^{n}x_{0},x_{0})-(l-\epsilon_{i})n.$
Since these numbers are unbounded, we can find a subsequence such
that $b_{i}(n_{i})>b_{i}(m)$ for any $m<n_{i}$ and by sequential
compactness we may moreover assume that $f^{n_{i}}(x_{0})\rightarrow P\in\overline{\mathcal{T}}^{GM}$,
cf. \cite{Ka01}. 

By a result of Liu and Su \cite{LS10} identifying the horoboundary
compactification of $(\mathcal{T},d)$ with the Gardiner-Masur compactification
(in particular showing that the latter is metrizable) we have for
any $k\geq1$ that \[
h_{P}(f^{k}x_{0})=\lim_{i\rightarrow\infty}d(f^{k}x_{0},f^{n_{i}}x_{0})-d(x_{0},f^{n_{i}}x_{0})\]

\[
\leq\liminf_{i\rightarrow\infty}d(x_{0},f^{n_{i}-k}x_{0})-d(x_{0},f^{n_{i}}x_{0})\]

\[
\leq\liminf_{i\rightarrow\infty}b_{i}(n_{i}-k)+(l-\epsilon_{i})(n_{i}-k)-b_{i}(n_{i})-(l-\epsilon_{i})n_{i}\]

\[
\leq\liminf_{i\rightarrow\infty}-(l-\epsilon_{i})k=-lk.\]
This means in terms of extremal lengths that\[
\left(\sup_{\beta}\frac{E_{P}(\beta)}{Ext_{f^{k}x_{0}}(\beta)^{1/2}}\right)^{-1}\geq Q(P)^{-1}e^{lk}.\]
and hence for any $\beta\in\mathcal{S}$ that \[
Ext_{f^{k}x_{0}}(\beta)\geq E_{P}(\beta)^{2}Q(P)^{-2}e^{2lk}.\]
On the other hand, in view of Kerckhoff's formula one has an estimate
from above:\[
e^{2d_{T}(f^{k}x_{0},x_{0})}=\sup_{\alpha}\frac{Ext_{f^{k}x_{0}}(\alpha)}{Ext_{x_{0}}(\alpha)}\geq\frac{Ext_{f^{k}x_{0}}(\beta)}{Ext_{x_{0}}(\beta)}.\]
In particular, provided $E_{P}(\beta)>0$, the two estimates imply
that \[
Ext_{f^{k}x_{0}}(\beta)^{1/n}\rightarrow e^{2l}\]
which by letting $\lambda=e^{2l}$ and setting $x=x_{0}$ finishes
the proof of Theorem \ref{thm:holomorphicNT}.

Note here that $l$, or $\lambda,$ are independent of $x$, while
$P$ might depend mildly on the point $x$.

\subsection*{Additional arguments concluding the proof of Theorem \ref{thm:DW}}

Suppose that the orbit is unbounded, then by a theorem of Calka, the
orbits tends to infinity, see \cite{Ka05}. Moreover, $l$ is defined
as above and we suppose that $P$ is a uniquely ergodic point. From
Miyachi \cite[Proposition 5.1]{Mi08} it then follows that we have
that this subsequence also converges to $P$ in the Thurston boundary.
Then we know from combining Theorem 11 with Corollary 45 in \cite{Ka05}
that the whole forward orbit must converge to $P$ (and hence in both
compactifications). 

In the notation above, we finally have\[
h_{P}(f(y))=\lim_{i\rightarrow\infty}d(f(y),f^{n_{i}}x_{0})-d(x_{0},f^{n_{i}}x_{0})\]
\[
\leq\liminf_{i\rightarrow\infty}d(y,f^{n_{i}-1}x_{0})-d(x_{0},f^{n_{i}}x_{0})\]
\[
\leq\liminf_{i\rightarrow\infty}d(y,f^{n_{i}-1}x_{0})-d(x_{0},f^{n_{i}-1}x_{0})-(l-\epsilon_{i})\]
\[
=h_{P}(y)-l\]
From this we may conclude Theorem \ref{thm:DW} as in the previous
proof (note that this time we have $\log Q(P)$ on both sides, hence
this cancels out).

In the bounded orbit case, we recall that if $f$ is a mapping class
then there is a fixed point by a theorem of Nielsen, if $f$ is a
strict contraction then by a simple argument due to Edelstein the
orbit converges towards the unique fixed point. However, in general
we do not know this.

We end by remarking that \cite[Proposition 5.1]{Mi08} also in the
reducible case (i.e. unbounded orbits but $P$ not uniquely ergodic)
gives information on the possible boundary limit points of the orbit
in terms of intersections with $P$. 

\noindent Anders Karlsson 

\noindent Section de mathématiques \\
Université de Genève

\noindent 2-4 Rue du Lièvre\\
Case Postale 64

\noindent 1211 Genève 4, 

\noindent Suisse 

\noindent e-mail: anders.karlsson@unige.ch 
\end{document}